\documentclass[10pt]{amsart}
\usepackage{mathtools,amssymb,bbm}
\usepackage{array,booktabs}
\usepackage{float}
\usepackage{tikz}
\usepackage[margin=1in]{geometry}
\usepackage{cite}

\usepackage{hyperref}
\usepackage{xcolor}

\hypersetup{
    colorlinks,
    linkcolor={red!50!black},
    citecolor={blue!50!black},
    urlcolor={blue!80!black},
    hypertexnames=false
}

\newtheorem{theorem}{Theorem}[section]
\newtheorem{lemma}[theorem]{Lemma}

\newtheorem{remark}[theorem]{Remark}

\newcommand{\tp}{{\scriptscriptstyle\mathsf{T}}}

\newcommand{\pp}{{\scriptscriptstyle++}}

\DeclareMathOperator{\vect}{vec}
\DeclareMathOperator{\svec}{svec}
\DeclareMathOperator{\mat}{mat}
\DeclareMathOperator{\smat}{smat}
\DeclareMathOperator{\Flag}{Flag}
\DeclareMathOperator{\V}{V}
\DeclareMathOperator{\tr}{tr}
\DeclareMathOperator{\Gr}{Gr}
\DeclareMathOperator{\im}{im}
\DeclareMathOperator{\spn}{span}
\DeclareMathOperator{\diag}{diag}

\DeclareMathOperator{\Cut}{\mathsf{cut}}
\DeclareMathOperator{\MaxCut}{\mathsf{maxcut}}
\DeclareMathOperator{\minimize}{minimize}

\newcommand{\ip}[2]{\left\langle #1,#2\right\rangle}

\newcommand{\half}{\frac12}

\newcommand{\E}{\mathcal E}
\newcommand{\A}{\mathcal A}

\let\S\undefined
\newcommand{\S}{\mathbb S}
\newcommand{\R}{\mathbb R}
\newcommand{\eps}{\varepsilon}
\newcommand{\one}{\mathbbm{1}}

\title{Quadratic programming on the Stiefel, Grassmann, and flag manifolds in polynomial time}
\author[T.~Tang]{Tianyun Tang}
\address{Computational and Applied Mathematics Initiative, Department of Statistics,
University of Chicago, Chicago, IL 60637}
\email{tianyuntang@uchicago.edu}
\author{Yuehaw Khoo}
\address{Computational and Applied Mathematics Initiative, Department of Statistics,
University of Chicago, Chicago, IL 60637}
\email{ykhoo@uchicago.edu}
\author[Z.~Lai]{Zehua~Lai}
\address{Institute of Applied Mathematics, Academy of Mathematics and Systems Science, Chinese Academy of Sciences, Beijing 100190, China}
\email{lzh@amss.ac.cn}
\author[L.-H.~Lim]{Lek-Heng~Lim}
\address{Computational and Applied Mathematics Initiative, Department of Statistics,
University of Chicago, Chicago, IL 60637}
\email{lekheng@uchicago.edu}

\begin{document}

\begin{abstract}
Quadratic programming over the Stiefel, Grassmannian, and flag manifolds is NP-hard in general, even though the corresponding linear optimization problems are easy.  We show that there is a broad and natural tractable region between these two facts.  If the quadratic part is relatively smaller than the spectral gap of the linear part, then the standard semidefinite relaxation will be exact: solving the SDP will solve the original nonconvex problem.  This gives polynomial-time algorithms for approximating the optimum value throughout this region.  We also prove converse hardness results showing that this tractable region is close to the best possible.  In short, the problem remains tractable when the separated linear term dominates, but NP-hard examples reappear once the quadratic term is large enough to overcome that separation.
\end{abstract}


\maketitle

\section{Introduction}

Optimization on matrix manifolds is now a standard part of continuous optimization.  Stiefel manifolds, Grassmannians, and flag manifolds occur naturally in dimension reduction, low-rank modeling, synchronization, signal processing, and related problems \cite{EdelmanAriasSmith1998,JourneeBachAbsilSepulchre2010,Singer2011,Manton2002}.  These manifolds are themselves nonconvex subsets of Euclidean space.  Thus even before a complicated objective is introduced, the feasible set already lies outside convex optimization.

The standard algorithms for manifold optimization, such as Riemannian gradient, Newton, and trust-region methods, are local in nature \cite{AbsilMahonySepulchre2008}.  They are highly effective in practice, but their usual guarantees concern stationarity or local optimality.  The present paper is about global optimization.  In this paper, we study when a nonconvex matrix-manifold problem can be solved to global optimality in polynomial time, and when this is hard in the usual complexity-theoretic sense.

We restrict attention to polynomial objectives.  This is a natural finite-dimensional class for matrix manifolds: the common models of \(\V(k,n)\), \(\Gr(k,n)\), and flag manifolds \(\Flag(k_1,\dots,k_p,n)\) are algebraic varieties, and polynomial objectives are specified by finitely many coefficients.  Among polynomial objectives, quadratic programming (QP)\footnote{In this paper, by QP we always mean quadratic programming with no constraints other than membership in the manifold.} is the first degree at which new global difficulty can appear.  Linear optimization on these manifolds is tractable, while Nesterov's work \cite{Nesterov2003} shows that cubic optimization on the unit sphere is already NP-hard.  It is therefore natural to ask what happens at degree two.  Recent work shows that QP is NP-hard in the worst case on Grassmannians and Stiefel manifolds \cite{LaiLimYe2025,LaiLimTang2026}.  Since Grassmannians are two-level flag manifolds, the Grassmannian result also implies worst-case NP-hardness for flag QP.

The three models studied in this paper are the following.  On the Stiefel manifold \(\V(k,n)=\{Y\in\R^{n\times k}:Y^\tp  Y=I_k\}\), we consider
\[
        \min_{Y\in\V(k,n)}
        \left\{\ip{A}{Y}+\half\vect(Y)^\tp  B\vect(Y)\right\}.
\]
On the Grassmannian, written as the set of rank-\(k\) orthogonal projectors $\Gr(k,n)=\{W\in\S^n:W^2=W,\ \tr W=k\}$, we consider
\[
        \min_{W\in\Gr(k,n)}
        \left\{\ip{A}{W}+\half\svec(W)^\tp  B\svec(W)\right\}.
\]
Finally, on a flag manifold \(\Flag(k_1,\dots,k_p,n)=\{U\Lambda U^\tp :U\in O(n)\}\), where \(0=k_0<k_1<\cdots<k_p<k_{p+1}=n\) and \(\Lambda=\diag(\lambda_1I_{k_1-k_0},\dots,\lambda_{p+1}I_{k_{p+1}-k_p})\) has the prescribed eigenvalue multiplicities, we consider
\[
        \min_{Z\in\Flag(k_1,\dots,k_p,n)}
        \left\{\ip{A}{Z}+\half\svec(Z)^\tp  B\svec(Z)\right\}.
\]
Here \(A\) denotes the linear coefficient in the appropriate ambient space: \(A\in\R^{n\times k}\) for Stiefel frames, and \(A\in\S^n\) for Grassmannians and flag manifolds.  Also, \(\vect\) is column-stacking vectorization, \(\svec\) is the scaled vectorization of symmetric matrices, and \(B\) denotes the Hessian matrix in the corresponding vectorization.

Worst-case hardness, however, is too coarse to explain the behavior of many structured instances.  The natural starting point is the linear case.  Over \(\V(k,n)\), a linear objective is solved by a singular value decomposition; over \(\Gr(k,n)\) and over \(\Flag(k_1,\dots,k_p,n)\), it is solved by an eigenvalue decomposition.  If the relevant singular value or eigenvalue spectral gap is positive, the optimizer of the linear problem is separated from its competitors.  This leads to a perturbative question:
\begin{quote}
If an easy linear problem is perturbed by a quadratic term that is small compared with this spectral separation, does the problem remain globally tractable, or can NP-hardness already appear in this regime?
\end{quote}

We prove that the separated linear term does control the problem below an explicit threshold.  For each manifold, we formulate a Shor semidefinite relaxation and prove that it is exact under a spectral gap condition.  The SDP has a unique rank-one optimizer, and its rank-one factor is the global minimizer of the original QP.  Consequently, the optimum value can be approximated to any prescribed additive accuracy in weak polynomial time; in exact arithmetic, solving the SDP also recovers a global optimizer.

We then show that this tractable region is close to the best possible.  At nearby scales beyond the sufficient conditions, one can construct instances for which an additive approximation to the optimum would recover the solution of an NP-hard problem.  For \(\Gr(k,n)\), the reduction uses bounded-degree stable-set instances together with the Motzkin--Straus formulation of the stable-set number \cite{MotzkinStraus1965,berman2002}.  For \(\V(k,n)\), the reduction uses a bounded-degree MaxCut instance embedded in diagonal coordinates \cite{berman2002}.  Since two-level flag manifolds are scaled copies of Grassmannians, the Grassmannian hardness result also gives hardness for flag QP.

Semidefinite relaxations, beginning with Shor's relaxation, are a standard approach to nonconvex QP \cite{Shor1987,VandenbergheBoyd1996}.  SDP relaxations for Stiefel QP have also been strengthened in recent work \cite{BurerPark2024}; that work focuses on relaxation strength and computation, rather than spectral-gap exactness.  Exactness of SDP relaxations has been studied from several viewpoints, including local stability of exact SDP relaxations \cite{CifuentesAgarwalParriloThomas2022}, exactness conditions for QCQPs \cite{WangKilincKarzan2022,JiangLi2023}, and exact recovery phenomena in synchronization problems \cite{BandeiraBoumalSinger2017}.  Our contribution is a deterministic, spectral gap-based exactness theorem for QP on three matrix manifolds, together with matching or nearly matching approximation-hardness results.

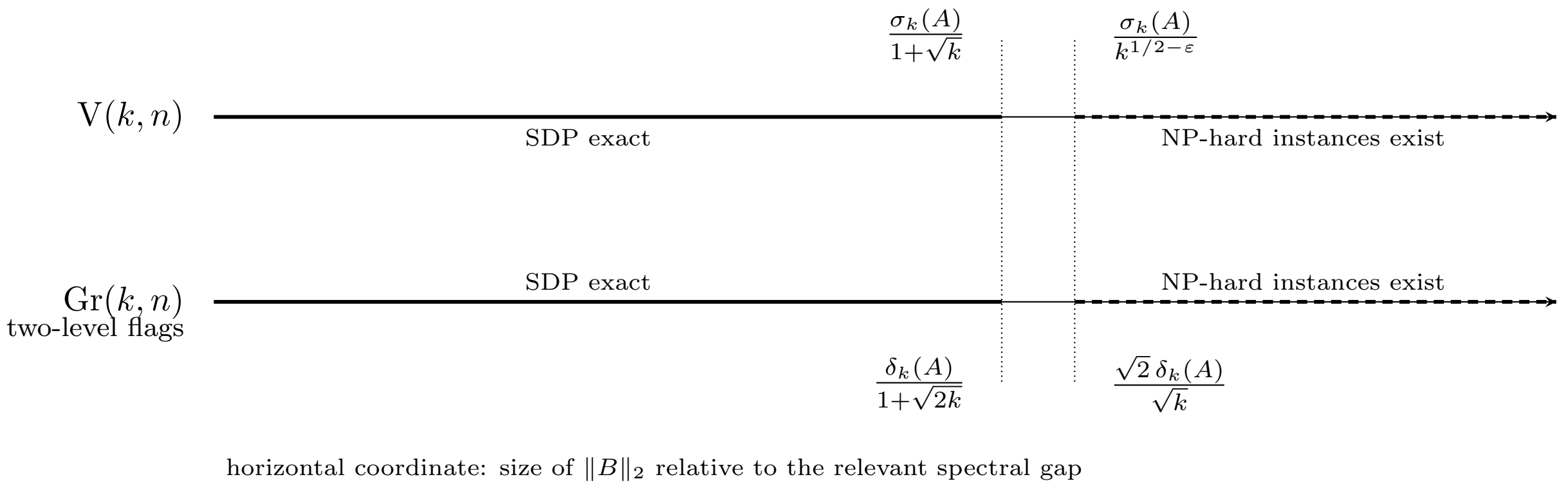
\begin{figure}[H]
\centering
\begin{tikzpicture}[x=1cm,y=1cm,>=stealth]
\def\xzero{0}
\def\xpoly{8.1}
\def\xhard{8.85}
\def\xend{13.8}
\def\yst{1.9}
\def\ygr{0}

\foreach \y in {\yst,\ygr} {
        \draw[->] (\xzero,\y) -- (\xend,\y);
        \draw[line width=1.1pt] (\xzero,\y) -- (\xpoly,\y);
        \draw[line width=1.1pt,densely dashed] (\xhard,\y) -- (\xend,\y);
}
\draw[densely dotted] (\xpoly,\ygr-0.82) -- (\xpoly,\yst+0.82);
\draw[densely dotted] (\xhard,\ygr-0.82) -- (\xhard,\yst+0.82);

\node[anchor=east] at (-0.18,\yst) {\(\V(k,n)\)};
\node[anchor=east] at (-0.18,\ygr) {\(\Gr(k,n)\)};
\node[anchor=east] at (-0.18,\ygr-0.28) {\footnotesize two-level flags};

\node[below] at (3.85,\yst) {\scriptsize SDP exact};
\node[below] at (11.2,\yst) {\scriptsize NP-hard instances exist};
\node[above] at (3.85,\ygr) {\scriptsize SDP exact};
\node[above] at (11.2,\ygr) {\scriptsize NP-hard instances exist};

\node[above,anchor=east] at (\xpoly-0.24,\yst+0.84) {\normalsize \(\frac{\sigma_k(A)}{1+\sqrt k}\)};
\node[above,anchor=west] at (\xhard+0.24,\yst+0.84) {\normalsize \(\frac{\sigma_k(A)}{k^{1/2-\eps}}\)};
\node[below,anchor=east] at (\xpoly-0.24,\ygr-0.84) {\normalsize \(\frac{\delta_k(A)}{1+\sqrt{2k}}\)};
\node[below,anchor=west] at (\xhard+0.24,\ygr-0.84) {\normalsize \(\frac{\sqrt2\,\delta_k(A)}{\sqrt k}\)};

\node[anchor=west] at (0,-1.72) {\scriptsize horizontal coordinate: size of \(\|B\|_2\) relative to the relevant spectral gap};
\end{tikzpicture}
\caption{Schematic comparison of the polynomial-time and NP-hardness scales.  Solid segments indicate regions where every instance is certified by the SDP relaxation.  Dashed segments indicate the scales at which the approximation-hardness reductions produce NP-hard families.}
\label{fig:regions}
\end{figure}

We use the following spectral-gap notation.  For the Stiefel linear coefficient \(A\in\R^{n\times k}\), \(\sigma_k(A)\) denotes the \(k\)-th singular value of \(A\), i.e., the smallest singular value in the compact singular value decomposition.  For a symmetric matrix \(A\in\S^n\), write its eigenvalues in increasing order,
\[
        \lambda_1(A)\le\dots\le\lambda_n(A),
\]
and define the Grassmannian gap
\[
        \delta_k(A)\coloneqq \lambda_{k+1}(A)-\lambda_k(A).
\]
For flag manifolds, we use the standard parameters
\(0=k_0<k_1<\cdots<k_p<k_{p+1}=n\) and
\[
        \Lambda=\diag(\lambda_1I_{k_1-k_0},\dots,
        \lambda_{p+1}I_{k_{p+1}-k_p}),
        \qquad \lambda_1>\cdots>\lambda_{p+1}.
\]
Thus the flag dimensions are the cumulative indices \(k_1,\dots,k_p\), while the eigenvalue multiplicities are the differences \(k_s-k_{s-1}\).  The relevant spectral cuts are
\[
        k_s,\qquad s=1,\dots,p,
\]
and the flag gap is 
\[
        \delta_*(A)\coloneqq \min_{1\le s\le p}
        \bigl(\lambda_{k_s+1}(A)-\lambda_{k_s}(A)\bigr).
\]
The positive and negative results are summarized in Table \ref{tab:main} and Figure \ref{fig:regions}.  The two columns should be read differently: the polynomial-time column gives a sufficient condition under which every instance is solved by the SDP relaxation, while the NP-hard column says that hard promised instances exist at the displayed scale.

\begin{table}[H]
\centering
\caption{Polynomial-time regions and hardness scales for QP on three manifolds.  The Stiefel hardness result holds for every fixed rational \(0<\eps<1/2\).  The Grassmannian hardness result holds for every fixed integer \(k\ge1\) and uses bounded-degree stable-set instances. The flag positive result applies to arbitrary flag manifolds \(\Flag(k_1,\dots,k_p,n)\), while flag hardness follows from the Grassmannian case.}
\setlength{\tabcolsep}{3pt}
\renewcommand{\arraystretch}{1.25}
\begin{tabular}{@{}>{\centering\arraybackslash}p{2.9cm}
                >{\centering\arraybackslash}p{3.1cm}
                >{\centering\arraybackslash}p{4.25cm}
                >{\centering\arraybackslash}p{4.3cm}@{}}
\toprule
manifold & spectral gap & polynomial-time & NP-hard \\
\midrule
\(\V(k,n)\) & \(\sigma_k(A)\) & \(\|B\|_2<\dfrac{\sigma_k(A)}{1+\sqrt k}\) &  \(\|B\|_2<\dfrac{\sigma_k(A)}{k^{1/2-\eps}}\) \\[3ex]
\(\Gr(k,n)\) & \(\delta_k(A)\) & \(\|B\|_2<\dfrac{\delta_k(A)}{1+\sqrt{2k}}\) &  \(\|B\|_2<\sqrt2\,\dfrac{\delta_k(A)}{\sqrt k}\) \\[3ex]
\(\Flag(k_1,\dots,k_p,n)\) & \(\delta_*(A)\) & \(\|B\|_2<\dfrac{\delta_*(A)}{C_\Lambda+\sqrt2\,E_\Lambda}\) & inherited from \(\Gr(k,n)\) \\
\bottomrule
\end{tabular}

\label{tab:main}
\end{table}

Here \(C_\Lambda\) and \(E_\Lambda\) depend only on the eigenvalue spacings of \(\Lambda\) and are defined in Section \ref{sec:flag}. For the standard Grassmannian projector normalization \(\Lambda=I_k\oplus0_{n-k}\), viewed as a two-level flag manifold, these constants reduce to \(C_\Lambda=1\) and \(E_\Lambda=\sqrt{k}\), so the flag exactness theorem generalizes the Grassmannian exactness theorem.

For flag manifolds, the positive theorem applies to arbitrary flag manifolds.  The hardness statement in the table is inherited from the two-level case: if \(\Lambda=\lambda_1I_k\oplus0_{n-k}\), then the two-level flag manifold \(\Flag(k,n)\) is a scaled copy of \(\Gr(k,n)\).  Thus Grassmannian hardness already rules out a general polynomial-time algorithm for flag QP.

The polynomial-time claim is meant in the weak SDP-optimization sense.  All input data are rational, and the input size is the total bit length of the rational entries and dimensions.  The requested accuracy is a positive rational number \(\eta\), encoded in binary.  A value algorithm returns a rational number within additive error \(\eta\) of the optimum in time polynomial in the input size and \(\log(1/\eta)\).  The relaxations below are explicit bounded spectrahedra.  Their affine equations and positive-semidefinite constraints admit standard rational weak separation oracles for affine slices of the positive semidefinite cone, and the weak separation/weak optimization equivalence then gives weak polynomial-time value optimization \cite{GroetschelLovaszSchrijver1988,VandenbergheBoyd1996}.  Thus the ellipsoid/cutting-plane framework solves them to additive accuracy \(\eta\) in time polynomial in the input size and \(\log(1/\eta)\).  

There are two notions here.  The exactness theorem is a mathematical statement about the relaxation: under the spectral gap conditions in Table \ref{tab:main}, the SDP optimizer is unique and rank one; in exact arithmetic, factoring this rank-one lift recovers the global optimizer of the original nonconvex problem.  The algorithmic consequence is slightly weaker, because standard SDP algorithms are approximate: they give an additive approximation to the optimum value in weak polynomial time.  The hardness results are stated already for value approximation.  Hence any algorithm that returns an additive approximate feasible point would imply a value approximation algorithm by evaluating the rational objective, to the precision specified in the reduction.

The paper is organized as follows.  Section \ref{sec:sdp} records the elementary SDP facts used later.  Sections \ref{sec:stiefel}, \ref{sec:grass}, and \ref{sec:flag} treat Stiefel, Grassmannian, and flag manifolds respectively.  The Stiefel and Grassmannian sections contain both polynomial-time and hardness theorems; the flag section gives the general polynomial-time theorem, with hardness inherited from Grassmannians as explained above.  Section \ref{sec:conclusion} concludes with the remaining endpoint question for Stiefel manifolds.

All matrices are real.  For matrices of the same size, \(\ip{X}{Y}=\tr(X^\tp  Y)\); for symmetric matrices this is \(\tr(XY)\).  The norm \(\|\cdot\|_2\) is the spectral/operator norm.  For \(X\in\R^{a\times b}\), \(\vect(X)\) denotes column-stacking vectorization.  For \(X\in\S^n\), \(\svec(X)\in\R^{n(n+1)/2}\) vectorizes the upper-triangular part of \(X\), with off-diagonal entries multiplied by \(\sqrt2\), so that
\[
        \ip{X}{Y}=\svec(X)^\tp \svec(Y).
\]
Its inverse is \(\smat\).

\section{Semidefinite relaxations and compactness}\label{sec:sdp}

The positive results have a common structure.  A positive spectral gap gives quadratic growth for the linear problem.  If the quadratic Hessian is sufficiently small, the global minimizer of the perturbed problem is still the optimizer of a nearby linear problem.  This nearby linear problem supplies the rank-one SDP certificate.  The three manifolds differ in the bookkeeping: singular values for Stiefel, eigenvalue gaps for Grassmannians, and several coupled eigenvalue gaps for flag manifolds.

We use the same elementary dual certificate in all three positive results.  The following lemma is the only linear-algebra fact needed for uniqueness.

\begin{lemma}[PSD range argument]\label{lem:psd-range}
Let \(S,H\succeq0\). If \(\ip{S}{H}=0\), then \(SH=0\), and hence \(\im(H)\subseteq\ker(S)\).  In particular, if \(\ker(S)=\spn\{u\}\) and the \((1,1)\) entry of \(H\) is one, then \(H=uu^\tp /u_1^2\), provided \(u_1\ne0\).
\end{lemma}

\begin{proof}
Since \(S^{1/2}HS^{1/2}\succeq0\) and has trace \(\ip{S}{H}=0\), it is zero.  Equivalently,
\[
        0=\tr(S^{1/2}HS^{1/2})=\|H^{1/2}S^{1/2}\|_F^2.
\]
Thus \(H^{1/2}S^{1/2}=0\), and hence \(SH=0\).  The last statement follows from the range inclusion and the normalization of the first entry.
\end{proof}

For the Grassmannian and flag relaxations, we use the following identities for scaled vectorization.  For each \(n\), put \(d=n(n+1)/2\), \(v_n\coloneqq \svec(I_n)\), and define \(\A_n:\S^d\to\R^d\) on symmetrized rank-one matrices by
\[
        \A_n\left(\frac{xz^\tp +zx^\tp }{2}\right)
        \coloneqq \half\svec(\smat(x)\smat(z)+\smat(z)\smat(x)),
\]
extended linearly. This definition is well-defined because the map \((x,z)\mapsto \half\svec(\smat(x)\smat(z)+\smat(z)\smat(x))\) is bilinear and symmetric, and hence induces a unique linear map on symmetric tensors. When the dimension is clear, we write \(\A\) and \(v\).

\begin{lemma}[symmetric lifting identities]\label{lem:svec-identities}
Let \(X,Y,Z\in\S^n\), and set \(x=\svec(X)\), \(y=\svec(Y)\), and \(z=\svec(Z)\).  Then
\[
        \A(xx^\tp )=\svec(X^2),
        \qquad
        \A^*(y)x=\svec\left(\frac{YX+XY}{2}\right),
\]
and
\[
        z^\tp \A^*(y)z=\ip{Y}{Z^2}.
\]
Moreover, for every \(X\in\S^d\), \(v^\tp \A(X)=\tr X\).
\end{lemma}

\begin{proof}
The first identity follows directly from the definition with \(z=x\).  For the adjoint identity, test against an arbitrary \(z=\svec(Z)\):
\[
\begin{aligned}
        z^\tp \A^*(y)x
        &=y^\tp \A\left(\frac{xz^\tp +zx^\tp }{2}\right)  \\
        &=\frac12\ip{Y}{XZ+ZX}
        =\ip{\frac{YX+XY}{2}}{Z}.
\end{aligned}
\]
This is the desired formula because \(\svec\) preserves the Frobenius inner product on symmetric matrices.  Taking \(X=Z\) gives \(z^\tp \A^*(y)z=\ip{Y}{Z^2}\).  Finally, \(v^\tp \A(xx^\tp )=\tr(X^2)=\|X\|_F^2=\tr(xx^\tp )\), and the identity extends to all \(X\in\S^d\) by linearity.
\end{proof}

The SDP feasible sets used below are compact.  We spell this out because it is the only regularity issue needed for polynomial-time solvability by SDP.

\begin{lemma}[compact SDP feasible sets]\label{lem:compact}
The Stiefel, Grassmannian, and flag Shor relaxations defined in Sections \ref{sec:stiefel}--\ref{sec:flag} have nonempty compact feasible sets.  Moreover their feasible matrices have trace bounded by a polynomial in the ambient dimension.
\end{lemma}

\begin{proof}
Nonemptiness follows from rank-one lifts of feasible manifold points.  For Stiefel, the constraint \(\E(X)=I_k\) gives \(\tr X=k\).  For Grassmannians, the equation \(\A(X)=x\), together with \(v^\tp x=k\), gives \(\tr X=k\) by Lemma \ref{lem:svec-identities}.  For flags, the same identity applied to each diagonal block \(X_{ss}\) gives \(\tr X_{ss}=k_s\), and positive semidefiniteness of the full lift bounds \(\tr X\) by the sum of the block traces.  A closed spectrahedron with bounded trace is compact.
\end{proof}

\section{The Stiefel manifold}\label{sec:stiefel}

Let
\[
        \V(k,n)=\{Y\in\R^{n\times k}:Y^\tp  Y=I_k\},\qquad 1\le k\le n.
\]
We consider
\begin{equation}\label{eq:stiefel-qp}
        \min_{Y\in\V(k,n)}
        \left\{\ip{A}{Y}+\half\vect(Y)^\tp  B\vect(Y)\right\},
\end{equation}
where \(A\in\R^{n\times k}\) and \(B\in\S^{nk}\).  Let \(\gamma\coloneqq \sigma_k(A)\).  The linear problem has a unique minimizer when \(\gamma>0\): if \(A=U\Sigma V^\tp \) is a compact singular value decomposition, then the minimizer is \(-UV^\tp \).

Partition \(X\in\S^{nk}\) into \(k\times k\) blocks \(X_{ij}\in\R^{n\times n}\), and define
\[
        \E(X)_{ij}=\tr(X_{ij}).
\]
The Shor relaxation is
\begin{equation}\label{eq:stiefel-sdp}
\begin{aligned}
\minimize_{y,X}\quad & \vect(A)^\tp  y+\half\ip{B}{X}\\
\text{subject to}\quad & \E(X)=I_k,\\
& \begin{bmatrix}1&y^\tp \\ y&X\end{bmatrix}\succeq0.
\end{aligned}\tag{S-St}
\end{equation}
This is the first-level lifting obtained by replacing \(yy^\tp \) with a free matrix \(X\) and linearizing the identities \(Y^\tp  Y=I_k\).  We first show that this elementary relaxation is already exact when the linear part has a sufficiently large singular-value spectral gap.

\begin{lemma}[linear growth on \(\V(k,n)\)]\label{lem:st-growth}
Let \(A=U\Sigma V^\tp \), \(\gamma=\sigma_k(A)>0\), and \(Y_A=-UV^\tp \).  Then every \(Y\in\V(k,n)\) satisfies
\[
        \ip{A}{Y-Y_A}\ge \frac{\gamma}{2}\|Y-Y_A\|_F^2.
\]
\end{lemma}

\begin{proof}
Set \(\widehat Y=YV\).  Since \(Y_AV=-U\),
\[
\ip{A}{Y-Y_A}=\tr\Sigma(U^\tp \widehat Y+I_k)
        =\sum_{i=1}^k \sigma_i(1+u_i^\tp  \widehat y_i).
\]
Each term \(1+u_i^\tp \widehat y_i\) is nonnegative.  Hence the expression is at least
\[
\gamma\sum_{i=1}^k(1+u_i^\tp \widehat y_i)=\frac\gamma2\|\widehat Y+U\|_F^2
        =\frac\gamma2\|Y-Y_A\|_F^2.
\]
\end{proof}

\begin{lemma}[Stiefel Shor certificate]\label{lem:st-certificate}
Let \(Y_*\in\V(k,n)\), \(y_*\coloneqq \vect(Y_*)\), and define
\[
        R_*\coloneqq A+\mat(By_*).
\]
Assume \(Y_*\) is the unique minimizer of \(Y\mapsto \ip{R_*}{Y}\) over \(\V(k,n)\).  Equivalently,
\[
        R_*=-Y_*P_*,\qquad P_*\in\S^k_\pp .
\]
If \(\lambda_{\min}(P_*)>\|B\|_2\), then the rank-one lift \((y_*,y_*y_*^\tp )\) is the unique optimizer of \eqref{eq:stiefel-sdp}.
\end{lemma}

\begin{proof}
The adjoint of \(\E\) is \(\E^*(P)=P\otimes I_n\).  The identity \(R_*=-Y_*P_*\) gives
\[
        \vect(A)+By_*+\E^*(P_*)y_*=0.
\]
Set
\[
        M\coloneqq \half(B+\E^*(P_*)),\qquad b\coloneqq \half\vect(A),\qquad \tau\coloneqq -y_*^\tp  M y_*.
\]
Then \(My_*+b=0\), and
\[
        M\succeq \frac{\lambda_{\min}(P_*)-\|B\|_2}{2}I\succ0.
\]
Define
\[
        S\coloneqq \begin{bmatrix}-\tau&b^\tp \\ b&M\end{bmatrix}
        =\begin{bmatrix}-y_*^\tp \\ I\end{bmatrix}M\begin{bmatrix}-y_*&I\end{bmatrix}\succeq0.
\]
We now verify the SDP KKT conditions.

\emph{Primal feasibility.}  The rank-one lift
\[
        H_*\coloneqq \begin{bmatrix}1&y_*^\tp \\ y_*&y_*y_*^\tp \end{bmatrix}
\]
is feasible for \eqref{eq:stiefel-sdp}, because \(\E(y_*y_*^\tp )=Y_*^\tp  Y_*=I_k\).

\emph{Dual feasibility.}  We identify \(S\) as the actual dual slack.  Introduce the lifted variable
\[
        H=\begin{bmatrix}t&y^\tp \\ y&X\end{bmatrix}
\]
and the affine constraint \(t=1\).  With the Lagrangian convention
\[
        L(H;\mu,\Lambda)=\vect(A)^\tp  y+\half\ip{B}{X}
        +\mu(t-1)+\ip{\Lambda}{\E(X)-I_k},
\]
choose \(\mu=-\tau\) and \(\Lambda=P_*/2\).  The coefficient of \(H\) in this Lagrangian is
\[
        \begin{bmatrix}
        \mu&b^\tp \\
        b&\half B+\E^*(\Lambda)
        \end{bmatrix}
        =
        \begin{bmatrix}
        -\tau&b^\tp \\
        b&M
        \end{bmatrix}
        =S.
\]
Thus the dual slack matrix is \(S\), and dual feasibility is exactly the already proved condition \(S\succeq0\).  The corresponding dual value is
\[
        -\mu-\ip{\Lambda}{I_k}
        =\tau-\half\tr P_* .
\]
Equivalently, for every feasible \(H\) of \eqref{eq:stiefel-sdp},
\[
        \vect(A)^\tp  y+\half\ip{B}{X}
        =\ip{S}{H}+\tau-\half\tr P_*
        \ge \tau-\half\tr P_* .
\]

\emph{Complementarity.}  Since \(My_*+b=0\), the factorization of \(S\) gives \(SH_*=0\), equivalently \(\ip{S}{H_*}=0\).  Therefore the primal value at \(H_*\) equals the dual value \(\tau-\half\tr P_*\), so \(H_*\) is optimal.

It remains only to prove uniqueness.  Since \(M\succ0\), \(\ker S=\spn\{(1,y_*^\tp )^\tp \}\).  If a feasible lift is optimal, then complementarity with \(S\) holds; Lemma \ref{lem:psd-range} forces its range to be contained in this one-dimensional kernel.  Hence the only optimal lift is \(H_*\).
\end{proof}

\begin{theorem}[SDP exactness for Stiefel QP]\label{thm:st-exact}
Assume \(\sigma_k(A)>0\).  If
\begin{equation}\label{eq:st-gap-bound}
        \|B\|_2<\frac{\sigma_k(A)}{1+\sqrt k},
\end{equation}
then the Shor relaxation \eqref{eq:stiefel-sdp} is exact and has a unique rank-one optimizer.  Consequently the optimum value of \eqref{eq:stiefel-qp} is computable by SDP to arbitrary prescribed additive accuracy in weak polynomial time.  In exact arithmetic, the global minimizer is recovered by factoring the unique rank-one SDP optimizer.
\end{theorem}

\begin{proof}
The argument is to convert a nonlinear global minimizer into the solution of a nearby linear Stiefel problem.  The singular-value spectral gap then gives the SDP certificate in Lemma \ref{lem:st-certificate}.

Put \(\gamma\coloneqq \sigma_k(A)\).
Let \(Y_*\) be a global minimizer of \eqref{eq:stiefel-qp}, and set
\[
        R_*\coloneqq A+\mat(B\vect(Y_*)).
\]
Since \(\|Y_*\|_F=\sqrt k\),
\[
        \sigma_k(R_*)\ge \gamma-\sqrt k\|B\|_2>\|B\|_2.
\]
Let \(\widehat Y\) be the unique minimizer of \(Y\mapsto\ip{R_*}{Y}\) over \(\V(k,n)\).  If \(\widehat Y\ne Y_*\), put \(\Delta=\widehat Y-Y_*\).  Lemma \ref{lem:st-growth}, applied to the linear matrix \(R_*\), gives
\[
        \ip{R_*}{\Delta}\le -\frac{\sigma_k(R_*)}{2}\|\Delta\|_F^2.
\]
Expanding the quadratic objective around \(Y_*\) gives
\[
        f(\widehat Y)-f(Y_*)=\ip{R_*}{\Delta}+\half\vect(\Delta)^\tp  B\vect(\Delta)
        \le -\frac{\sigma_k(R_*)-\|B\|_2}{2}\|\Delta\|_F^2<0,
\]
contradicting the optimality of \(Y_*\).  Therefore \(Y_*\) is the polar minimizer of \(R_*\).  Lemma \ref{lem:st-certificate} applies with \(\lambda_{\min}(P_*)=\sigma_k(R_*)>\|B\|_2\).
\end{proof}

We now turn to the NP-hardness region for Stiefel QP.  In the hardness statements, an additive \(\eta\)-approximate solution means a feasible point whose objective value is within \(\eta\) of the global optimum; for maximization this means value at least \(\operatorname{OPT}-\eta\), and for minimization it means value at most \(\operatorname{OPT}+\eta\).

The Stiefel hardness construction does not reach the endpoint \(\sigma_k(A)/\sqrt k\), but it gets arbitrarily close in exponent.

\begin{theorem}[NP-hardness for Stiefel QP]\label{thm:st-hard}
Fix a rational number \(0<\eps<1/2\).  It is NP-hard to approximate the minimum value within additive error \(1\) for Stiefel QP over the family of manifolds \(\V(k,n)\), \(n\ge k\), with \(k,n\) part of the input, even restricted to rational instances satisfying
\[
        \|B\|_2<\frac{\sigma_k(A)}{k^{1/2-\eps}}.
\]
Consequently, computing an additive \(1\)-approximate solution is NP-hard as well.
\end{theorem}

The theorem is a variable-dimension hardness statement.  In the reduction below, \(k\) is chosen as a polynomial function of the graph size.  This is enough to prove NP-hardness for the general Stiefel QP family, because an algorithm that solved all rational instances with arbitrary input dimensions \(k,n\) would in particular solve this subfamily.  We do not claim here that the same hardness result holds for every fixed \(k\).

\begin{proof}
It is convenient to write the reduction as a maximization problem first.  Multiplying the resulting objective by \(-1\) gives the minimization instance in the theorem and leaves \(\|B\|_2\) and \(\sigma_k(A)\) unchanged.

The reduction uses diagonal entries of the frame.  The first \(s\) diagonal variables form a padding block with a large positive linear coefficient.  The last \(m\) diagonal variables encode a MaxCut instance.  A bilinear cancellation term removes the unwanted linear contribution on the MaxCut block once the padding variables are forced to equal one.

We reduce from MaxCut on graphs of maximum degree at most three, which is NP-hard by the hardness of \(3\)-MAX CUT \cite[Theorem~1(i)]{berman2002}.  Let \(G\) have \(m\) vertices, adjacency matrix \(A_G\), and edge set \(E(G)\).  Put
\[
        C_G\coloneqq I_m-A_G.
\]
Since \(G\) has maximum degree at most three, all eigenvalues of the symmetric matrix \(A_G\) lie in \([-3,3]\), so the eigenvalues of \(C_G=I_m-A_G\) lie in \([-2,4]\).  Hence \(\|C_G\|_2\le4\).  Moreover, for \(z\in\{\pm1\}^m\),
\begin{equation}\label{eq:maxcut-identity}
        z^\tp  C_Gz=m-2|E(G)|+4\Cut_G(z).
\end{equation}

Choose an integer \(q\) with \(2\eps q>1\), and an integer \(M\ge2\) with \(M^{2\eps}\ge 32\).  Set
\[
        k\coloneqq Mm^q,\qquad s\coloneqq k-m.
\]
The role of \(\eps\) is to make the padding dimension \(k\) large enough compared with the graph size \(m\).  Since \(2\eps q>1\) and \(M^{2\eps}\ge32\), we have that
\begin{equation}\label{eq:st-padding}
        m\le \frac{1}{32}k^{2\eps}.
\end{equation}
 Also \(s\ge k/2\) and \(s\ge m\). Choose an integer \(\gamma\) such that
\begin{equation}\label{eq:gamma-choice}
        \gamma\ge 32k^{1/2-\eps}.
\end{equation}
Let \(n\ge k\) be arbitrary, and work over \(\V(k,n)\).  For \(Y\in\V(k,n)\), define
\[
        u_j\coloneqq Y_{jj}\ (1\le j\le s),\qquad z_i\coloneqq Y_{s+i,s+i}\ (1\le i\le m),
\]
\[
        S\coloneqq \sum_{j=1}^s u_j,
        \qquad T\coloneqq \sum_{i=1}^m z_i.
\]
Let
\[
        A=\begin{bmatrix}
        \diag(\smash[t]{\overbrace{2\gamma,\dots,2\gamma}^{s\text{ times}},
        \overbrace{\gamma,\dots,\gamma}^{m\text{ times}})}\\
        0_{(n-k)\times k}
        \end{bmatrix},
\]
so \(\sigma_k(A)=\gamma\), and define the quadratic term by
\begin{equation}\label{eq:st-hard-obj}
        \half\vect(Y)^\tp  B\vect(Y)=z^\tp  C_Gz-\frac{\gamma}{s}ST.
\end{equation}
The matrix \(B\) is zero outside the coordinates \(Y_{jj}\), \(1\le j\le k\).
The objective for maximization is
\[
        \Psi(u,z)=2\gamma S+\gamma T+z^\tp  C_Gz-\frac{\gamma}{s}ST,
        \qquad (u,z)\in[-1,1]^s\times[-1,1]^m.
\]
Indeed, every \(Y\in\V(k,n)\) has \(|Y_{jj}|\le1\), so the Stiefel objective is bounded above by this box maximization problem.
For fixed \(z\), the coefficient of every \(u_j\) is
\[
        2\gamma-\frac{\gamma}{s}T\ge 2\gamma-\frac{\gamma m}{s}\ge\gamma>0.
\]
Hence every box maximizer has \(u_1=\dots=u_s=1\).  On this face the linear term in \(T\) cancels, and the box problem reduces to
\[
        2\gamma s+\max_{z\in[-1,1]^m} z^\tp  C_Gz.
\]
The latter maximum is attained at a sign vector: with all other variables fixed, the objective is a convex quadratic function of each coordinate because \((C_G)_{ii}=1\).  Every such sign vector is realized in \(\V(k,n)\) by the embedded diagonal frame
\[
        Y=\begin{bmatrix}
        \diag(\smash[t]{\overbrace{1,\dots,1}^{s\text{ times}},z_1,\dots,z_m)}\\
        0_{(n-k)\times k}
        \end{bmatrix}.
\]
Thus the box upper bound is attained by a feasible Stiefel point.  Therefore, by \eqref{eq:maxcut-identity}, the exact Stiefel optimum is
\[
        2\gamma s+m-2|E(G)|+4\MaxCut(G),
\]
and determines \(\MaxCut(G)\).
In fact, an additive \(1\)-approximation already determines \(\MaxCut(G)\): if \(C\coloneqq 2\gamma s+m-2|E(G)|\), then the possible optimum values are \(C+4r\), \(r\in\mathbb Z\), and these values are separated by \(4\).  Given an approximate value \(\widehat V\) with \(|\widehat V-\operatorname{OPT}|\le1\), the integer
\[
        \operatorname{round}\left(\frac{\widehat V-C}{4}\right)
\]
is exactly \(\MaxCut(G)\).  If an algorithm returns an additive \(1\)-approximate maximizer instead, evaluating its rational objective value to error less than \(1\) gives a total value error strictly less than \(2\).  After division by \(4\), the error is strictly less than \(1/2\), so the same rounding recovers \(\MaxCut(G)\).  Thus an additive \(1\)-approximate maximizer would also recover \(\MaxCut(G)\).

It remains to check the promise.  The active quadratic matrix, in variables \((u,z)\), is
\[
        \begin{bmatrix}
        0_s&-(\gamma/s)\one_s\one_m^\tp \\
        -(\gamma/s)\one_m\one_s^\tp &2C_G
        \end{bmatrix}.
\]
The full matrix \(B\) is obtained by embedding this active matrix into the diagonal-coordinate subspace and putting zeros elsewhere.  Hence its norm is the norm of the displayed matrix.  By the triangle inequality,
\[
        \|B\|_2
        \le
        \left\|
        \begin{bmatrix}
        0_s&0\\
        0&2C_G
        \end{bmatrix}
        \right\|_2
        +
        \left\|
        \begin{bmatrix}
        0_s&-(\gamma/s)\one_s\one_m^\tp \\
        -(\gamma/s)\one_m\one_s^\tp &0_m
        \end{bmatrix}
        \right\|_2.
\]
The first norm is \(2\|C_G\|_2\).  The second norm is
\[
        \frac{\gamma}{s}\|\one_s\one_m^\tp \|_2
        =\gamma\sqrt{\frac{m}{s}}.
\]
Therefore
\[
        \|B\|_2\le 2\|C_G\|_2+\gamma\sqrt{\frac{m}{s}}
        \le 8+\gamma\sqrt{\frac{m}{s}}.
\]
By \eqref{eq:st-padding},
\[
        \sqrt{\frac{m}{s}}\le \sqrt{\frac{2m}{k}}\le \frac{1}{4}k^{-1/2+\eps},
\]
and by \eqref{eq:gamma-choice},
\[
        8\le \frac{1}{4}\frac{\gamma}{k^{1/2-\eps}}.
\]
Therefore
\[
        \|B\|_2\le \frac{1}{2}\frac{\gamma}{k^{1/2-\eps}}<\frac{\sigma_k(A)}{k^{1/2-\eps}}.
\]
All parameters have polynomial bit length for fixed \(\eps\).  Since MaxCut remains NP-hard on graphs of maximum degree at most three \cite[Theorem~1(i)]{berman2002}, the promised Stiefel approximation problem is NP-hard.
\end{proof}

\begin{remark}
The proof shows why the endpoint \(\sigma_k(A)/\sqrt k\) is not reached.  If the hard MaxCut block uses \(m\) diagonal variables and the padding block has size \(s\), the cancellation term has norm proportional to \(\sigma_k(A)\sqrt{m/s}\).  Since \(m\) must grow in this reduction, one loses a factor \(m^{1/2}\) relative to the endpoint.
\end{remark}

\section{The Grassmannian}\label{sec:grass}

Let
\[
        \Gr(k,n)=\{W\in\S^n:W^2=W,\ \tr W=k\}.
\]
We consider
\begin{equation}\label{eq:grass-qp}
        \min_{W\in\Gr(k,n)}
        \left\{\ip{A}{W}+\half\svec(W)^\tp  B\svec(W)\right\}.
\end{equation}
Let the eigenvalues of \(A\) be ordered increasingly, and put
\[
        \delta_k(A)\coloneqq \lambda_{k+1}(A)-\lambda_k(A).
\]
When \(\delta_k(A)>0\), the linear problem has the unique minimizer given by the projector onto the \(k\) smallest eigenvalues of \(A\).

Let \(v\coloneqq \svec(I_n)\), and let \(\A\) be the map in Lemma \ref{lem:svec-identities} for this dimension.
The Shor relaxation is
\begin{equation}\label{eq:grass-sdp}
\begin{aligned}
\minimize_{x,X}\quad & \svec(A)^\tp  x+\half\ip{B}{X}\\
\text{subject to}\quad & \A(X)=x,\\
& v^\tp  x=k,\\
& Xv=kx,\\
& \begin{bmatrix}1&x^\tp \\x&X\end{bmatrix}\succeq0.
\end{aligned}\tag{S-Gr}
\end{equation}
The equation \(Xv=kx\) is valid for every rank-one lift: if \(X=xx^\tp \), then \(Xv=x(x^\tp  v)=kx\).  This inexpensive additional equation strengthens the basic Shor relaxation and helps keep it aligned with the rank-one geometry of projectors.

\begin{lemma}[projector growth]\label{lem:grass-growth}
Let \(W_A\) be the rank-\(k\) projector onto the \(k\) smallest eigenvalues of \(A\), and assume \(\delta_k(A)>0\).  Then every \(W\in\Gr(k,n)\) satisfies
\[
        \ip{A}{W-W_A}\ge \frac{\delta_k(A)}{2}\|W-W_A\|_F^2.
\]
\end{lemma}

\begin{proof}
Choose an eigenbasis of \(A\) adapted to \(W_A\).  We may assume that \(A\) is diagonal and that the first \(k\) diagonal entries of \(W_A\) are equal to one.  Then
\begin{align*}
\ip{A}{W-W_A}&=\sum_{i=1}^k A_{ii}(W_{ii}-1) + \sum_{i=k+1}^n A_{ii}W_{ii}\\
&\ge \lambda_k(A)\sum_{i=1}^k (W_{ii}-1) + \lambda_{k+1}(A) \sum_{i=k+1}^n W_{ii} \\
&= (\lambda_{k+1}(A)-\lambda_k(A)) \sum_{i=k+1}^n W_{ii} \\
&= \delta_k(A)\tr((I-W_A)W)\\
&=\frac{\delta_k(A)}{2}\|W-W_A\|_F^2. \qedhere
\end{align*}
\end{proof}

\begin{lemma}[Grassmannian Shor certificate]\label{lem:grass-certificate}
Let \(W_*\in\Gr(k,n)\), \(x_*\coloneqq \svec(W_*)\), and set
\[
        A_*\coloneqq A+\smat(Bx_*).
\]
Suppose \(W_*\) is the minimum rank-\(k\) spectral projector of \(A_*\), with spectral gap
\[
        \delta_*\coloneqq \lambda_{k+1}(A_*)-\lambda_k(A_*).
\]
If \(\delta_*>\|B\|_2\), then the rank-one lift of \(W_*\) is the unique optimizer of \eqref{eq:grass-sdp}.
\end{lemma}

\begin{proof}
Write \(A_*=U\diag(\eta_1,\dots,\eta_n)U^\tp \), \(\eta_1\le\dots\le\eta_n\), with
\[
        W_*=U\diag(I_k,0)U^\tp .
\]
Let \(\alpha=(\eta_k+
\eta_{k+1})/2\), and define
\[
        R=U\diag(\eta_1-\alpha,
        \dots,
        \eta_k-\alpha,
        -\eta_{k+1}+\alpha,
        \dots,
        -\eta_n+
        \alpha)U^\tp .
\]
Then
\begin{equation}\label{eq:grass-separation}
        A_*-(RW_*+W_*R-R)-\alpha I=0,
        \qquad -R\succeq \frac{\delta_*}{2}I.
\end{equation}
Let \(r\coloneqq \svec(R)\).  Applying \(\svec\) to \eqref{eq:grass-separation} gives
\[
        \svec(A)+Bx_*-2\A^*(r)x_*+r-\alpha v=0.
\]
Set
\[
        M\coloneqq \half B-\A^*(r),\qquad
        b\coloneqq \half(\svec(A)+r-\alpha v),
        \qquad \tau\coloneqq -x_*^\tp  Mx_*.
\]
Then \(Mx_*+b=0\).  Moreover, by Lemma \ref{lem:svec-identities}, for \(z=\svec(Z)\),
\[
        z^\tp [-\A^*(r)]z=-\ip{R}{Z^2}\ge \frac{\delta_*}{2}\|Z\|_F^2,
\]
so
\[
        M\succeq \frac{\delta_*-
\|B\|_2}{2}I\succ0.
\]
Define
\[
        S\coloneqq \begin{bmatrix}-\tau&b^\tp \\ b&M\end{bmatrix}
        =\begin{bmatrix}-x_*^\tp \\ I\end{bmatrix}M\begin{bmatrix}-x_*&I\end{bmatrix}\succeq0.
\]
We verify the SDP KKT conditions.

\emph{Primal feasibility.}  The rank-one lift
\[
        H_*\coloneqq \begin{bmatrix}1&x_*^\tp \\ x_*&x_*x_*^\tp \end{bmatrix}
\]
is feasible for \eqref{eq:grass-sdp}: Lemma \ref{lem:svec-identities} gives \(\A(x_*x_*^\tp )=x_*\), while \(v^\tp  x_*=k\) and \(x_*x_*^\tp  v=kx_*\).

\emph{Dual feasibility.}  We again check that \(S\) has the dual-slack form.  Write
\[
        H=\begin{bmatrix}t&x^\tp \\ x&X\end{bmatrix}
\]
and include the affine constraint \(t=1\).  Write the Grassmannian constraints as
\[
        x-\A(X)=0,\qquad k-v^\tp  x=0,\qquad Xv-kx=0.
\]
Using multipliers \(\mu=-\tau\), \(r\), \(\alpha\), and \(0\), respectively, the Lagrangian is
\[
\begin{aligned}
        L(H;\mu,r,\alpha)
        &=\svec(A)^\tp  x+\half\ip{B}{X}
        +\mu(t-1)\\
        &\quad+\ip{r}{x-\A(X)}
        +\alpha(k-v^\tp  x).
\end{aligned}
\]
The coefficient of \(H\) in this Lagrangian is
\[
        \begin{bmatrix}
        \mu&\half(\svec(A)+r-\alpha v)^\tp \\
        \half(\svec(A)+r-\alpha v)&\half B-\A^*(r)
        \end{bmatrix}
        =
        \begin{bmatrix}
        -\tau&b^\tp \\
        b&M
        \end{bmatrix}
        =S.
\]
Thus \(S\) is the dual slack matrix, and dual feasibility is the condition \(S\succeq0\).  The corresponding dual value is
\[
        -\mu+\alpha k=\tau+\alpha k.
\]
Equivalently, for every feasible \(H\) of \eqref{eq:grass-sdp},
\[
        \svec(A)^\tp  x+\half\ip{B}{X}
        =\ip{S}{H}+\tau+\alpha k
        \ge \tau+\alpha k.
\]

\emph{Complementarity.}  Since \(Mx_*+b=0\), the factorization of \(S\) gives \(SH_*=0\), equivalently \(\ip{S}{H_*}=0\).  Hence \(H_*\) attains the dual value and is optimal.

Since \(M\succ0\), \(\ker S=\spn\{(1,x_*^\tp )^\tp \}\).  Any optimal feasible lift must be complementary to \(S\), so Lemma \ref{lem:psd-range} gives uniqueness.
\end{proof}

\begin{theorem}[SDP exactness for Grassmannian QP]\label{thm:grass-exact}
Assume \(\lambda_{k+1}(A)>\lambda_k(A)\).  If
\begin{equation}\label{eq:grass-gap-bound}
        \|B\|_2<\frac{\lambda_{k+1}(A)-\lambda_k(A)}{1+\sqrt{2k}},
\end{equation}
then the Shor relaxation \eqref{eq:grass-sdp} is exact and has a unique rank-one optimizer.  Consequently the optimum value of \eqref{eq:grass-qp} is computable by SDP to arbitrary prescribed additive accuracy in weak polynomial time.  In exact arithmetic, the global minimizer is recovered by factoring the unique rank-one SDP optimizer.
\end{theorem}

\begin{proof}
As in the Stiefel case, we absorb the derivative of the quadratic term into the linear matrix and then use the spectral gap of the resulting linear problem.

Put \(\delta\coloneqq \lambda_{k+1}(A)-\lambda_k(A)\).
Let \(W_*\) be a global minimizer and set \(A_*\coloneqq A+\smat(B\svec(W_*))\). Write \(E=A_*-A\).  For a symmetric perturbation, Weyl's inequality gives
\[
        \lambda_{k+1}(A_*)-\lambda_k(A_*)
        \ge \delta-\bigl(\lambda_{\max}(E)-\lambda_{\min}(E)\bigr).
\]
Moreover \(\lambda_{\max}(E)-\lambda_{\min}(E)\le \sqrt2\|E\|_F\), since \((a-b)^2\le2(a^2+b^2)\) for the largest and smallest eigenvalues.  As \(\svec\) is an isometry and \(\|W_*\|_F=\sqrt k\),
\[
        \|E\|_F=\|\smat(B\svec(W_*))\|_F
        \le \|B\|_2\sqrt k.
\]
Therefore
\[
        \lambda_{k+1}(A_*)-\lambda_k(A_*)
        \ge \delta-\sqrt{2k}\|B\|_2>\|B\|_2.
\]
Let \(\widehat W\) be the minimum rank-\(k\) spectral projector of \(A_*\).  If \(\widehat W\ne W_*\), put \(\Delta=\widehat W-W_*\).  Expanding around \(W_*\), and using Lemma \ref{lem:grass-growth} for \(A_*\),
\[
\begin{aligned}
        f(\widehat W)-f(W_*)
        &=\ip{A_*}{\Delta}+\half\svec(\Delta)^\tp  B\svec(\Delta)\\
        &\le -\frac{\lambda_{k+1}(A_*)-\lambda_k(A_*)-
        \|B\|_2}{2}\|\Delta\|_F^2<0,
\end{aligned}
\]
contradicting optimality.  Hence \(W_*\) is the minimum spectral projector of \(A_*\), and Lemma \ref{lem:grass-certificate} applies.
\end{proof}

We now prove that the scale \((\lambda_{k+1}(A)-\lambda_k(A))/\sqrt k\) is unavoidable in a worst-case complexity sense.  The proof uses only diagonal entries of the projector.  By Schur--Horn, the possible diagonals of rank-\(k\) projectors form the hypersimplex, and the Motzkin--Straus theorem then lets us place a stable-set problem inside this diagonal slice.

\begin{lemma}[diagonal image and stable sets]\label{lem:diag-stable}
The diagonal image of \(\Gr(k,n)\) is
\[
        \{w\in[0,1]^n:\one^\tp  w=k\}.
\]
If \(G\) is a graph with adjacency matrix \(A_G\), \(C_G\coloneqq I+A_G\), and independence number \(\alpha(G)\), then
\[
        \min_{y\in\Delta_m}y^\tp  C_Gy=\frac1{\alpha(G)},
        \qquad
        \min_{x\ge0,\ \one^\tp  x=r}x^\tp  C_Gx=\frac{r^2}{\alpha(G)}.
\]
\end{lemma}

\begin{proof}
The diagonal statement is the Schur-Horn theorem applied to the eigenvalue vector \((1,\dots,1,0,\dots,0)\).  The stable-set formula is Motzkin-Straus \cite{MotzkinStraus1965} applied to the complement graph.  For \(y\in\Delta_m\),
\[
        y^\tp (I+A_G)y
        =1-2\sum_{\{i,j\}\in E(\overline G)}y_iy_j,
\]
and \(\omega(\overline G)=\alpha(G)\).  The scaled formula follows by writing \(x=ry\).
\end{proof}

\begin{theorem}[NP-hardness for Grassmannian QP]\label{thm:grass-hard}
For every fixed integer \(k\ge1\), there is a polynomial-time reduction from maximum independent set red problem on graphs of red $m$ vertices and maximum degree at most three to rational instances of \eqref{eq:grass-qp} on \(\Gr(k,k+m)\) satisfying
\begin{equation}\label{eq:grass-hard-bound}
        \|B\|_2<\sqrt2\,\frac{\delta_k(A)}{\sqrt k}.
\end{equation}
For these instances, it is NP-hard to compute an additive \(1\)-approximate solution, or even to approximate the optimum value within additive error \(1\). 
\end{theorem}

\begin{proof}
We reduce from maximum independent set on graphs of maximum degree at most three, which is NP-hard by the inapproximability of \(3\)-MIS \cite[Theorem~1(v)]{berman2002}.  Let \(G\) be such a graph on \(m\) vertices, with adjacency matrix \(A_G\), independence number \(a_G=\alpha(G)\), and
\[
        C_G\coloneqq I+A_G.
\]
Since \(G\) has maximum degree at most three,
\[
        \|C_G\|_2\le \|C_G\|_\infty\le4.
\]
We first treat the case \(k\ge2\).  Set
\[
        L\coloneqq 392m(m+1).
\]
The coefficient \(3/4\) in the \(\tau^2\) term below is used only to give a strict margin in the \(k=2\) norm estimate; it also keeps the one-dimensional scalar minimizer inside \([0,1]\).  Because \(k\) is fixed, choose an integer \(q=q_k>1\), depending only on \(k\), such that
\begin{equation}\label{eq:q-choice}
        \max\left\{8,\ \frac34+\sqrt{\frac9{16}+\frac{(q-1)^2}{k-1}}\right\}<\frac{\sqrt2\,q}{\sqrt k}.
\end{equation}
Such a \(q\) exists.  If \(k>2\), then the second term on the left has leading coefficient \(1/\sqrt{k-1}\), while the right-hand side has leading coefficient \(\sqrt{2/k}\).  Since \(1/\sqrt{k-1}<\sqrt{2/k}\) exactly when \(k>2\), the inequality holds for all sufficiently large \(q\).  If \(k=2\), then \eqref{eq:q-choice} becomes
\[
        \max\left\{8,\ \frac34+\sqrt{\frac9{16}+(q-1)^2}\right\}<q.
\]
The first inequality is \(8<q\).  For the second inequality, since \(q>1\), we may square after moving \(3/4\) to the right:
\[
\begin{aligned}
        \frac34+\sqrt{\frac9{16}+(q-1)^2}<q
        &\iff \frac9{16}+(q-1)^2<\left(q-\frac34\right)^2\\
        &\iff q>2.
\end{aligned}
\]
Thus, for example, any integer \(q>8\) works when \(k=2\).  Now set
\[
        \delta\coloneqq qL,\qquad b\coloneqq \frac{\delta-L}{k-1}.
\]
For \(W\in\Gr(k,k+m)\), write
\[
        u_j=W_{jj}\ (1\le j\le k-1),\qquad
        \tau=W_{kk},\qquad
        x_i=W_{k+i,k+i}\ (1\le i\le m),
\]
\[
        U\coloneqq \sum_{j=1}^{k-1}u_j.
\]
By Lemma \ref{lem:diag-stable}, the diagonal variables satisfy
\begin{equation}\label{eq:grass-hard-polytope}
        0\le u_j,\tau,x_i\le1,
        \qquad U+\tau+\one^\tp x=k.
\end{equation}
Let \(A\) be diagonal with entries \(-4\delta\) on coordinates \(1,\dots,k-1\), entry \(-\delta\) on coordinate \(k\), and entries \(0\) on coordinates \(k+1,\dots,k+m\).  Then \(\delta_k(A)=\delta\).  Define the quadratic term by
\begin{equation}\label{eq:grass-hard-quadratic}
        \half\svec(W)^\tp B\svec(W)=Lx^\tp C_Gx+\frac34 L\tau^2+bU\tau.
\end{equation}
The objective over the diagonal polytope is
\[
        F(u,\tau,x)=-4\delta U-\delta\tau+Lx^\tp C_Gx+\frac34 L\tau^2+bU\tau.
\]
Every minimizer has \(U=k-1\).  Indeed, if \(U<k-1\), then some \(u_j<1\), and \eqref{eq:grass-hard-polytope} gives \(\tau+\one^\tp x>1\).  If \(\tau>0\), moving a small amount from \(\tau\) to such a coordinate \(u_j\) has directional derivative
\[
        -4\delta+\delta-\frac32 L\tau+b(\tau-U)\le -3\delta+b<0.
\]
If \(\tau=0\), then some \(x_i>0\).  Moving a small amount from this \(x_i\) to \(u_j\) gives derivative
\[
        -4\delta+b\tau-2L(C_Gx)_i=-4\delta-2L(C_Gx)_i<0,
\]
because \(C_G\ge0\) and \(x\ge0\).  Hence no minimizer can have \(U<k-1\).

On the face \(U=k-1\), the remaining variables satisfy \(\tau+\one^\tp x=1\), and $-\delta \tau+bU\tau=-L\tau$ by the definition of the constants. Therefore, by Lemma \ref{lem:diag-stable}, the problem becomes
\[
        -4\delta(k-1)+L\min_{0\le\tau\le1}
        \left\{\frac{(1-\tau)^2}{a_G}-\tau+\frac34\tau^2\right\}.
\]
For \(s=1/a_G\), define
\[
        g(s)\coloneqq \min_{0\le\tau\le1}\{s(1-\tau)^2-\tau+\frac34\tau^2\}.
\]
Differentiating the quadratic in \(\tau\), the minimizer is interior and equals \(\tau_s=(4s+2)/(4s+3)\).  Substitution gives the closed form \(g(s)=-(s+1)/(4s+3)\).  By the envelope theorem, or by direct differentiation after substituting \(\tau_s\),
\[
        g'(s)=(1-\tau_s)^2=\frac{1}{(4s+3)^2}\ge\frac1{49},\qquad 0<s\le1.
\]
Thus for \(1\le t\le m-1\),
\[
        g\!\left(\frac1t\right)-g\!\left(\frac1{t+1}\right)
        \ge \frac1{49}\left(\frac1t-\frac1{t+1}\right)
        \ge \frac1{49m(m+1)}.
\]
The exact optimum is
\begin{equation}\label{eq:grass-hard-value}
        \min_{W\in\Gr(k,k+m)}F(W)=-R(k-1)+Lg(1/\alpha(G)).
\end{equation}
Since \(L=392m(m+1)\), the possible values in \eqref{eq:grass-hard-value} are separated by at least \(8\).  Hence an additive \(1\)-approximation to the optimum value determines \(\alpha(G)\) by rounding to the nearest number in the list \(-R(k-1)+Lg(1/t)\), \(t=1,\dots,m\).  If an algorithm returns an additive \(1\)-approximate minimizer instead, evaluating its rational objective value to error less than \(1\) still leaves total error less than \(2\), below half of the separation, so the same rounding recovers \(\alpha(G)\).

It remains to check the promise.  In the diagonal variables, the active quadratic matrix is
\[
        Q=
        \begin{bmatrix}
        0_{k-1}&b\one_{k-1}&0\\
        b\one_{k-1}^\tp &\frac32 L&0\\
        0&0&2LC_G
        \end{bmatrix}.
\]
The \(x\)-block has norm at most \(8L\).  The \((u,\tau)\)-block has eigenvalues \(0\), with multiplicity \(k-2\), and
\[
       \frac34 L\pm\sqrt{\frac9{16}L^2+b^2(k-1)}.
\]
Therefore, using \(b=(q-1)L/(k-1)\) and \eqref{eq:q-choice},
\[
\begin{aligned}
        \|B\|_2
        &\le L\max\left\{8,\ \frac34+\sqrt{\frac9{16}+\frac{(q-1)^2}{k-1}}\right\}\\
        &< \frac{\sqrt2\,qL}{\sqrt k}
        =\sqrt2\,\frac{\delta_k(A)}{\sqrt k}.
\end{aligned}
\]

It remains to handle \(k=1\).  Work on \(\Gr(1,m+1)\), write \(\tau=W_{11}\) and \(x_i=W_{1+i,1+i}\).  Use the same \(L=392m(m+1)\), set \(\delta=6L\), \(M=4L\), and take \(A=-\delta e_1e_1^\tp\), so \(\delta_1(A)=\delta\).  Minimize
\[
        Lx^\tp C_Gx-\delta\tau+M\tau^2.
\]
For fixed \(\tau\), this equals \(L\{s(1-\tau)^2-6\tau+4\tau^2\}\), with \(s=1/\alpha(G)\).  Its minimizer is \(\tau_s=(s+3)/(s+4)\), and the closed-form value function is \(g_1(s)=-(2s+9)/(s+4)\).  Moreover,
\[
        g_1'(s)=(1-\tau_s)^2=\frac1{(s+4)^2}\ge\frac1{25},\qquad 0<s\le1.
\]
Thus the possible optimum values are separated by at least \(L/(25m(m+1))>8\), so additive error \(1\) again determines \(\alpha(G)\).  Finally,
\[
        \|B\|_2\le \max\{2M,2L\|C_G\|_2\}\le8L<6\sqrt2 L
        =\sqrt2\,\delta_1(A).
\]
This proves the promised bound also in rank one, and completes the reduction.
\end{proof}

\section{Flag manifolds}\label{sec:flag}

We use the standard flag notation from the matrix-orbit model.  Fix integers
\[
        0=k_0<k_1<\cdots<k_p<k_{p+1}=n
\]
and distinct real numbers \(\lambda_1>\cdots>\lambda_{p+1}\).  Put
\[
        \Lambda=\diag(\lambda_1 I_{k_1-k_0},
        \lambda_2 I_{k_2-k_1},\dots,
        \lambda_{p+1} I_{k_{p+1}-k_p})\in\S^n.
\]
Thus the parameters \(k_1,\dots,k_p\) are the dimensions of the nested subspaces, while the multiplicities of the eigenvalue blocks of \(\Lambda\) are the differences \(k_s-k_{s-1}\).  The flag manifold in this embedded model is
\[
        \Flag(k_1,\dots,k_p,n)=\{U\Lambda U^\tp :U\in O(n)\}.
\]
We assume without loss of generality that the smallest eigenvalue is zero.  Indeed, if the smallest eigenvalue of the original orbit matrix is \(c\), write \(Z=cI+\widetilde Z\).  Then the objective in \(Z\) becomes a quadratic objective in \(\widetilde Z\), with linear coefficient \(A+\smat(B\svec(cI))\) and with only an additive constant changed.  Thus replacing \(\Lambda\) by \(\Lambda-cI\) preserves the class of QP instances and gives \(\lambda_{p+1}=0\).
This orbit display describes the flag manifold but is not the coordinate system used for the QP below.  If one substitutes \(Z=U\Lambda U^\tp\) directly into a quadratic objective in \(Z\), the entries of \(Z\) are already quadratic in \(U\), so the resulting optimization problem is quartic in the orthogonal variable \(U\).  The first-level Shor relaxation used in this paper applies instead to a quadratic objective in embedded matrix coordinates.

Put \(\rho_s=\lambda_s-\lambda_{s+1}>0\), \(s=1,\dots,p\).  Every \(Z\in\Flag(k_1,\dots,k_p,n)\) can be written as
\begin{equation}\label{eq:flag-projectors}
        Z=\Phi(D),
        \qquad
        \Phi(D)\coloneqq \sum_{s=1}^{p}\rho_sD_s,
\end{equation}
where
\[
        D_s^2=D_s,
        \qquad \tr D_s=k_s,
        \qquad D_sD_{s+1}=D_s,\quad s=1,\dots,p-1.
\]
Conversely, every nested tuple \((D_1,\dots,D_p)\) gives a point of \(\Flag(k_1,\dots,k_p,n)\).  This is the product of Grassmannians embedding of the flag manifold used in \cite[Proposition~3]{YeWongLim2022}.  In these variables the map \(D\mapsto Z\) is affine, so a quadratic objective in \(Z\) remains a quadratic objective in the projector variables \(D_s\).  We study
\begin{equation}\label{eq:flag-qp}
        \min_{Z\in\Flag(k_1,\dots,k_p,n)}
        \left\{\ip{A}{Z}+\half\svec(Z)^\tp  B\svec(Z)\right\}.
\end{equation}
In the projector variables this is
\begin{equation}\label{eq:flag-proj-qp}
        \min_D
        \left\{\sum_{s=1}^{p}\rho_s\ip{A}{D_s}
        +\half\svec(\Phi(D))^\tp  B\svec(\Phi(D))\right\}.
\end{equation}
Define
\[
        \delta_s(A)\coloneqq \lambda_{k_s+1}(A)-\lambda_{k_s}(A),
        \qquad
        \delta_*(A)\coloneqq \min_{1\le s\le p}\delta_s(A),
\]
\[
        C_\Lambda\coloneqq \sum_{s=1}^{p}\rho_s,
        \qquad
        E_\Lambda\coloneqq \sum_{s=1}^{p}\rho_s\sqrt{k_s}.
\]

Let \(d=n(n+1)/2\) and \(v=\svec(I_n)\).  The Shor relaxation uses variables \(x_s=\svec(D_s)\in\R^d\), \(x=(x_1,\dots,x_p)\), and a lifted matrix \(X=(X_{st})_{s,t=1}^{p}\in\S^{pd}\).  Define the multiplication map \(\mathcal N:\R^{d\times d}\to\R^{n^2}\) and the linear map \(L:\R^d\to\R^{n^2}\) by
\[
        \mathcal N(ab^\tp )\coloneqq \vect(\smat(a)\smat(b)),
        \qquad
        L(a)\coloneqq \vect(\smat(a)),
\]
extended linearly.  The function $\mathcal{N}$ is well-defined because the right-hand side is bilinear in \((a,b)\), and hence induces a unique linear map on matrices.
For a rank-one lift, \(X_{s,s+1}=x_sx_{s+1}^\tp\).  Therefore the product \(D_sD_{s+1}=D_s\) linearizes to
\[
        \mathcal N(X_{s,s+1})=L(x_s),
        \qquad s=1,\dots,p-1.
\]
The flag Shor relaxation is
\begin{equation}\label{eq:flag-sdp}
\begin{aligned}
\minimize_{x,X}\quad & \sum_s \rho_s\svec(A)^\tp  x_s
        +\half\ip{\widehat B}{X}\\
\text{subject to}\quad & \A(X_{ss})=x_s,
        \quad v^\tp  x_s=k_s,
        \quad X_{ss}v=k_sx_s,\qquad s=1,\dots,p,\\
& \mathcal N(X_{s,s+1})=L(x_s),
        \qquad s=1,\dots,p-1,\\
& \begin{bmatrix}1&x^\tp \\x&X\end{bmatrix}\succeq0.
\end{aligned}\tag{S-Fl}
\end{equation}
Here \(\widehat B\) is the block Hessian induced by \(\svec(\Phi(D))^\tp  B\svec(\Phi(D))\), namely \(\widehat B_{st}=\rho_s\rho_t B\). The positive theorem below is therefore a multi-cut spectral-gap argument.

In the proof we also use the product problem obtained by dropping the nesting constraints from \eqref{eq:flag-proj-qp}:
\begin{equation}\label{eq:flag-product-qp}
\begin{aligned}
\minimize_{D_1,\dots,D_p}\quad &
        F(D)
        \coloneqq \sum_{s=1}^{p}\rho_s\ip{A}{D_s}
        +\half\svec(\Phi(D))^\tp  B\svec(\Phi(D))\\
\text{subject to}\quad &
        D_s\in\Gr(k_s,n),\qquad s=1,\dots,p .
\end{aligned}\tag{P-Fl}
\end{equation}
The corresponding product Shor relaxation keeps the same block lift and cross-block objective, but drops the linearized nesting equations:
\begin{equation}\label{eq:flag-product-sdp}
\begin{aligned}
\minimize_{x,X}\quad & \sum_s \rho_s\svec(A)^\tp  x_s
        +\half\ip{\widehat B}{X}\\
\text{subject to}\quad & \A(X_{ss})=x_s,
        \quad v^\tp  x_s=k_s,
        \quad X_{ss}v=k_sx_s,
        \qquad s=1,\dots,p,\\
& \begin{bmatrix}1&x^\tp \\x&X\end{bmatrix}\succeq0.
\end{aligned}\tag{S-P-Fl}
\end{equation}
The flag SDP relaxation \eqref{eq:flag-sdp} is a strengthening of \eqref{eq:flag-product-sdp}.

\begin{lemma}[projector growth]\label{lem:proj-growth}
Let \(D_C\) be the rank-\(N\) projector onto the \(N\) smallest eigenvalues of \(C\in\S^n\), and assume \(\lambda_{N+1}(C)>\lambda_N(C)\).  Then every \(D\in\Gr(N,n)\) satisfies
\[
        \ip{C}{D-D_C}\ge
        \frac{\lambda_{N+1}(C)-\lambda_N(C)}{2}
        \|D-D_C\|_F^2.
\]
\end{lemma}

\begin{proof}
Choose an eigenbasis of \(C\) adapted to \(D_C\).  We may assume that \(C\) is diagonal, with eigenvalues ordered increasingly, and that
\[
       D_C=\diag(I_N,0).
\]
Since \(D\) is an orthogonal projector, \(0\preceq D\preceq I\).  Hence \(0\le D_{ii}\le1\) for every \(i\), so \(D_{ii}-1\le0\) for \(i\le N\) and \(D_{ii}\ge0\) for \(i>N\).  Therefore
\begin{align*}
\ip{C}{D-D_C}
&=\sum_{i=1}^N \lambda_i(C)(D_{ii}-1)+\sum_{i=N+1}^n \lambda_i(C)D_{ii}\\
&\ge \lambda_N(C)\sum_{i=1}^N(D_{ii}-1)
        +\lambda_{N+1}(C)\sum_{i=N+1}^nD_{ii}\\
&=(\lambda_{N+1}(C)-\lambda_N(C))\sum_{i=N+1}^nD_{ii}\\
&=(\lambda_{N+1}(C)-\lambda_N(C))\tr((I-D_C)D)\\
&=\frac{\lambda_{N+1}(C)-\lambda_N(C)}{2}\|D-D_C\|_F^2.
\end{align*}
Here the third line uses \(\tr D=N\), and the last line uses that both \(D\) and \(D_C\) are rank-\(N\) orthogonal projectors.
\end{proof}

\begin{lemma}[weighted flag growth]\label{lem:flag-growth}
Let \(C\in\S^n\) have positive gaps at all cuts \(k_s\), and let \(D_s^C\) be its minimum rank-\(k_s\) spectral projectors.  Then for any \(\widetilde D_s\in\Gr(k_s,n)\),
\[
        \sum_{s=1}^{p}\rho_s\ip{C}{\widetilde D_s-D_s^C}
        \ge \frac{\delta_*(C)}{2}
        \sum_{s=1}^{p}\rho_s\|\widetilde D_s-D_s^C\|_F^2.
\]
\end{lemma}

\begin{proof}
Apply Lemma \ref{lem:proj-growth} to each cut and multiply by \(\rho_s\).
\end{proof}

\begin{theorem}[SDP exactness for flag QP]\label{thm:flag-exact}
Assume \(\delta_*\coloneqq \delta_*(A)>0\).  If
\begin{equation}\label{eq:flag-gap-bound}
        \|B\|_2<\frac{\delta_*}{C_\Lambda+\sqrt2\,E_\Lambda},
\end{equation}
then the flag Shor relaxation \eqref{eq:flag-sdp} is exact and has a unique rank-one optimizer.  Consequently the optimum value of \eqref{eq:flag-qp} is computable by SDP to arbitrary prescribed additive accuracy in weak polynomial time.  In exact arithmetic, the global minimizer is recovered by factoring the unique rank-one SDP optimizer.
\end{theorem}

\begin{proof}
We divide the proof into three steps.  The first step shows that the product relaxation has an automatically nested minimizer.  The second step constructs the spectral certificates explicitly. The third step assembles these certificates into the Shor dual certificate.

\emph{Step 1: reduction to the product problem.}
It is enough to certify the product Shor relaxation \eqref{eq:flag-product-sdp}; the flag SDP relaxation \eqref{eq:flag-sdp} is a strengthening and contains the same rank-one flag point.

Let \(D_*=(D_{1,*},\dots,D_{p,*})\) be a global minimizer of the product problem \eqref{eq:flag-product-qp} and put
\begin{equation}\label{eq:flag-proof-Astar}
        A_*\coloneqq A+\smat\bigl(B\svec(\Phi(D_*))\bigr).
\end{equation}
Write \(E=A_*-A\).  As in the Grassmannian proof,
\begin{equation}\label{eq:flag-proof-gap-Astar}
\begin{aligned}
        \delta_*(A_*)
        &\ge \delta_*-\bigl(\lambda_{\max}(E)-\lambda_{\min}(E)\bigr)\\
        &\ge \delta_*-\sqrt2\|E\|_F\\
        &\ge \delta_*-\sqrt2\,E_\Lambda\|B\|_2
        > C_\Lambda\|B\|_2 .
\end{aligned}
\end{equation}
The third inequality uses \(\|E\|_F\le\|B\|_2\|\Phi(D_*)\|_F\le E_\Lambda\|B\|_2\).
For each \(s\), let \(\widehat D_s\) be the minimum rank-\(k_s\) spectral projector of \(A_*\), and set
\begin{equation}\label{eq:flag-proof-Delta}
        \Delta_s\coloneqq \widehat D_s-D_{s,*}.
\end{equation}
Expanding \(F\) at \(D_*\) gives
\begin{equation}\label{eq:flag-proof-expansion}
\begin{aligned}
F(\widehat D)-F(D_*)
&=\sum_s\rho_s\ip{A_*}{\Delta_s}
+\frac12\svec\left(\sum_s\rho_s\Delta_s\right)^\tp B
\svec\left(\sum_s\rho_s\Delta_s\right).
\end{aligned}
\end{equation}
The linear term in \eqref{eq:flag-proof-expansion} is controlled by Lemma~\ref{lem:flag-growth}, applied to the common matrix \(A_*\):
\begin{equation}\label{eq:flag-proof-linear-control}
        \sum_s\rho_s\ip{A_*}{\Delta_s}
        \le -\frac{\delta_*(A_*)}{2}
        \sum_s\rho_s\|\Delta_s\|_F^2.
\end{equation}
The quadratic term is controlled by the weighted Cauchy--Schwarz inequality
\begin{equation}\label{eq:flag-proof-weighted-cauchy}
        \left\|\sum_s\rho_s\Delta_s\right\|_F^2
        \le \left(\sum_s\rho_s\right)
        \left(\sum_s\rho_s\|\Delta_s\|_F^2\right)
        = C_\Lambda\sum_s\rho_s\|\Delta_s\|_F^2.
\end{equation}
Combining \eqref{eq:flag-proof-expansion}, \eqref{eq:flag-proof-linear-control}, and \eqref{eq:flag-proof-weighted-cauchy}, we obtain
\begin{equation}\label{eq:flag-proof-negative-drop}
\begin{aligned}
F(\widehat D)-F(D_*)
&\le -\frac{\delta_*(A_*)-C_\Lambda\|B\|_2}{2}
        \sum_s\rho_s\|\Delta_s\|_F^2.
\end{aligned}
\end{equation}
By \eqref{eq:flag-proof-gap-Astar}, the coefficient in \eqref{eq:flag-proof-negative-drop} is strictly negative unless all \(\Delta_s=0\).  Hence every \(D_{s,*}\) is the minimum rank-\(k_s\) spectral projector of the same matrix \(A_*\).  These projectors are therefore nested, so the product minimizer is an exact flag point.

\emph{Step 2: explicit construction of the spectral certificates.}
Choose an eigenvalue decomposition
\begin{equation}\label{eq:flag-proof-eigendecomp}
        A_*=U\diag(\eta_1,\dots,\eta_n)U^\tp,
        \qquad \eta_1\le\cdots\le\eta_n,
\end{equation}
with the basis chosen so that
\begin{equation}\label{eq:flag-proof-Ds-diagonal}
        D_{s,*}=U\diag(I_{k_s},0)U^\tp,
        \qquad s=1,\dots,p.
\end{equation}
This is possible because the projectors in Step~1 are the nested spectral projectors of the same matrix \(A_*\).  For each \(s\), define
\begin{equation}\label{eq:flag-proof-alpha}
        \alpha_s\coloneqq \frac{\eta_{k_s}+\eta_{k_s+1}}{2}
\end{equation}
and
\begin{equation}\label{eq:flag-proof-Rs-def}
        R_s\coloneqq U\diag(r_{s,1},\dots,r_{s,n})U^\tp,
        \qquad
        r_{s,i}=\begin{cases}
        \eta_i-\alpha_s, & i\le k_s,\\
        \alpha_s-\eta_i, & i>k_s.
        \end{cases}
\end{equation}
Then, on the first \(k_s\) coordinates, \(R_sD_{s,*}+D_{s,*}R_s-R_s\) has diagonal entries \(r_{s,i}=\eta_i-
\alpha_s\); on the remaining coordinates it has diagonal entries \(-r_{s,i}=\eta_i-\alpha_s\).  Therefore
\begin{equation}\label{eq:flag-proof-cut-certificate}
        A_*-(R_sD_{s,*}+D_{s,*}R_s-R_s)-\alpha_sI=0.
\end{equation}
Moreover, for \(i\le k_s\),
\[
        -r_{s,i}=\alpha_s-\eta_i\ge \alpha_s-\eta_{k_s}
        =\frac{\eta_{k_s+1}-\eta_{k_s}}{2},
\]
and for \(i>k_s\),
\[
        -r_{s,i}=\eta_i-\alpha_s\ge \eta_{k_s+1}-\alpha_s
        =\frac{\eta_{k_s+1}-\eta_{k_s}}{2}.
\]
Thus
\begin{equation}\label{eq:flag-proof-Rs-bound}
        -R_s\succeq \frac{\delta_s(A_*)}{2}I
        \succeq \frac{\delta_*(A_*)}{2}I.
\end{equation}

\emph{Step 3: the product Shor certificate.}
Let \(q_s\coloneqq \svec(R_s)\), \(x_{s,*}\coloneqq \svec(D_{s,*})\), and \(x_*=(x_{1,*},\dots,x_{p,*})\).  Define
\begin{equation}\label{eq:flag-proof-M-def}
        M\coloneqq \frac12\widehat B-
        \diag\bigl(\A^*(\rho_1q_1),\dots,\A^*(\rho_pq_p)\bigr),
\end{equation}
\begin{equation}\label{eq:flag-proof-btau-def}
        b\coloneqq \frac12\bigl(\rho_1(\svec(A)+q_1-\alpha_1v),\dots,
        \rho_p(\svec(A)+q_p-\alpha_pv)\bigr),
        \qquad
        \tau\coloneqq -x_*^\tp Mx_* .
\end{equation}
The identities \eqref{eq:flag-proof-cut-certificate} imply
\begin{equation}\label{eq:flag-proof-stationarity}
        Mx_*+b=0.
\end{equation}
Indeed, the \(s\)-th block of \(\widehat Bx_*\) is
\begin{equation}\label{eq:flag-proof-Bhat-block}
        (\widehat Bx_*)_s
        =\rho_sB\svec(\Phi(D_*)),
\end{equation}
and the \(s\)-th block of \eqref{eq:flag-proof-stationarity} is one half of the \(\svec\)-form of
\begin{equation}\label{eq:flag-proof-block-stationarity}
        \rho_s\{A_*-(R_sD_{s,*}+D_{s,*}R_s-R_s)-\alpha_sI\}=0.
\end{equation}
Here we used Lemma~\ref{lem:svec-identities}, which gives
\(2\A^*(q_s)x_{s,*}=\svec(R_sD_{s,*}+D_{s,*}R_s)\).

Next we prove that \(M\succ0\).  For \(z=(z_1,\dots,z_p)\), set \(Z_s=\smat(z_s)\) and \(Z_\rho=\sum_s\rho_sZ_s\).  By Lemma~\ref{lem:svec-identities} and \eqref{eq:flag-proof-Rs-bound},
\begin{equation}\label{eq:flag-proof-M-positive}
\begin{aligned}
        z^\tp Mz
        &=\frac12\svec(Z_\rho)^\tp B\svec(Z_\rho)
        -\sum_s\rho_s\ip{R_s}{Z_s^2}  \\
        &\ge -\frac12\|B\|_2\|Z_\rho\|_F^2
        +\frac{\delta_*(A_*)}{2}\sum_s\rho_s\|Z_s\|_F^2 \\
        &\ge \frac{\delta_*(A_*)-C_\Lambda\|B\|_2}{2}
        \sum_s\rho_s\|z_s\|_2^2>0,
\end{aligned}
\end{equation}
where the last line again uses \eqref{eq:flag-proof-weighted-cauchy}.  Because all \(\rho_s\) are positive, the weighted norm \(\sum_s\rho_s\|z_s\|_2^2\) is positive for every nonzero \(z\).  Thus \(M\succ0\).

Now set
\begin{equation}\label{eq:flag-proof-slack}
        S\coloneqq \begin{bmatrix}-\tau&b^\tp \\ b&M\end{bmatrix}
        =\begin{bmatrix}-x_*^\tp \\ I\end{bmatrix}
        M\begin{bmatrix}-x_*&I\end{bmatrix}\succeq0.
\end{equation}
We verify the KKT conditions for the product Shor relaxation \eqref{eq:flag-product-sdp}.

\emph{Primal feasibility.}  The rank-one lift
\begin{equation}\label{eq:flag-proof-rankone-lift}
        H_*\coloneqq \begin{bmatrix}1&x_*^\tp \\ x_*&x_*x_*^\tp \end{bmatrix}
\end{equation}
is feasible for \eqref{eq:flag-product-sdp}, because each \(D_{s,*}\) is a rank-\(k_s\) projector.  Since the projectors are nested, the same lift is also feasible for the stronger flag relaxation \eqref{eq:flag-sdp}.

\emph{Dual feasibility.}  Let
\begin{equation}\label{eq:flag-proof-H-variable}
        H=\begin{bmatrix}t&x^\tp \\ x&X\end{bmatrix}
\end{equation}
and include the affine constraint \(t=1\).  For each \(s\), write the product-SDP constraints as
\begin{equation}\label{eq:flag-proof-product-constraints}
        x_s-\A(X_{ss})=0,
        \qquad k_s-v^\tp x_s=0,
        \qquad X_{ss}v-k_sx_s=0.
\end{equation}
Use multipliers \(\mu=-\tau\), \(\rho_sq_s\), \(\rho_s\alpha_s\), and \(0\), respectively.  The Lagrangian is
\begin{equation}\label{eq:flag-proof-lagrangian}
\begin{aligned}
        L(H)
        &=\sum_s\rho_s\svec(A)^\tp x_s
        +\frac12\ip{\widehat B}{X}
        +\mu(t-1) \\
        &\quad+\sum_s\ip{\rho_sq_s}{x_s-\A(X_{ss})}
        +\sum_s\rho_s\alpha_s(k_s-v^\tp x_s).
\end{aligned}
\end{equation}
By the definitions \eqref{eq:flag-proof-M-def} and \eqref{eq:flag-proof-btau-def}, the coefficient matrix of \(H\) in \eqref{eq:flag-proof-lagrangian} is exactly the matrix \(S\) in \eqref{eq:flag-proof-slack}.  Therefore \(S\succeq0\) is a dual-feasible slack matrix, and every feasible point of \eqref{eq:flag-product-sdp} satisfies
\begin{equation}\label{eq:flag-proof-dual-identity}
        \sum_s\rho_s\svec(A)^\tp x_s
        +\frac12\ip{\widehat B}{X}
        =\ip{S}{H}+\tau+\sum_s\rho_s\alpha_sk_s
        \ge \tau+\sum_s\rho_s\alpha_sk_s.
\end{equation}

\emph{Complementarity and uniqueness.}  Since \(Mx_*+b=0\), the factorization in \eqref{eq:flag-proof-slack} gives
\begin{equation}\label{eq:flag-proof-complementarity}
        SH_*=0,
        \qquad \text{equivalently}\qquad \ip{S}{H_*}=0.
\end{equation}
Thus \(H_*\) attains the dual value in \eqref{eq:flag-proof-dual-identity} and is optimal for the product Shor relaxation.  Since \(M\succ0\),
\begin{equation}\label{eq:flag-proof-kernel}
        \ker S=\spn\{(1,x_*^\tp)^\tp\}.
\end{equation}
Any optimal feasible lift for the product SDP is complementary to \(S\).  Lemma~\ref{lem:psd-range} then forces its range to be contained in the one-dimensional space in \eqref{eq:flag-proof-kernel}; because its \((1,1)\) entry is one, it must equal \(H_*\).  Hence the product Shor relaxation has the unique optimizer \(H_*\).

Finally, \eqref{eq:flag-sdp} is a strengthening of \eqref{eq:flag-product-sdp} and contains \(H_*\).  Therefore its optimum value is the same as the product SDP optimum, and every optimal point of \eqref{eq:flag-sdp} is also product-SDP optimal.  By the uniqueness just proved, the flag Shor relaxation also has the unique rank-one optimizer \(H_*\).  This proves exactness and the stated polynomial-time value consequence follows from Lemma~\ref{lem:compact} and the SDP weak-optimization discussion in the introduction.

\end{proof}

Finally, flag QP is NP-hard in the worst case because Grassmannians are two-level flag manifolds.  Indeed, if
\[
        \Lambda=\lambda_1I_k\oplus0_{n-k},\qquad \lambda_1>0,
\]
then every \(Z\in\Flag(k,n)\) has the form
\[
        Z=\lambda_1W,\qquad W\in\Gr(k,n).
\]
Given a Grassmannian QP
\[
        g(W)=\ip{A_G}{W}+\half\svec(W)^\tp  B_G\svec(W),
\]
define a flag objective on this two-level orbit by
\[
        B_F\coloneqq \lambda_1^{-2}B_G,\qquad
        \svec(A_F)\coloneqq \lambda_1^{-1}\svec(A_G).
\]
Then substituting \(Z=\lambda_1W\) gives
\[
        \ip{A_F}{Z}+\half\svec(Z)^\tp  B_F\svec(Z)=g(W).
\]
Thus an additive \(\eta\)-approximation for the flag instance gives an additive \(\eta\)-approximation for the Grassmannian instance.  The Grassmannian approximation-hardness theorem therefore already implies NP-hardness for approximating flag QP in the worst case.

\section{Conclusion}\label{sec:conclusion}

QP on the Stiefel manifold, Grassmannian, and flag manifolds is NP-hard in the worst case.  Nevertheless, a separated linear term creates a polynomial-time region.  The reason is simple: the spectral gap gives the linear problem enough stability along the manifold, and a quadratic perturbation below the spectral gap-scaled threshold cannot destroy the rank-one SDP certificate.

For Grassmannians and flag manifolds, the polynomial-time and NP-hardness bounds differ only by constant factors: after the sharpened exactness estimate, the two constants on the \(\delta_k(A)/\sqrt k\) scale are asymptotically \(1/\sqrt2\) and \(\sqrt2\).  The hardness proof uses bounded-degree stable-set instances and many fixed directions to cancel a large linear field at spectral-norm cost \(\delta_k(A)/\sqrt k\).  For Stiefel manifolds, we prove NP-hardness at \(\sigma_k(A)/k^{1/2-\eps}\) for every fixed rational \(\eps>0\), but this does not reach \(\sigma_k(A)/\sqrt k\).  Determining the exact Stiefel endpoint remains an interesting open problem.

\bibliographystyle{amsplain}
\providecommand{\bysame}{\leavevmode\hbox to3em{\hrulefill}\thinspace}
\providecommand{\MR}{\relax\ifhmode\unskip\space\fi MR }
\providecommand{\MRhref}[2]{%
  \href{http://www.ams.org/mathscinet-getitem?mr=#1}{#2}
}
\providecommand{\href}[2]{#2}

\end{document}